\documentclass [12pt]{article}
\usepackage{amsfonts}
\usepackage{amssymb}
\usepackage{multirow}
\usepackage{amsmath,amssymb,amsthm}
\newtheorem{theorem}{Theorem}

\newtheorem{proposition}{Proposition}
\newtheorem{corollary}{Corollary}
\newcommand{\re}[1]{\mbox{${\rm (\ref{#1})}$}}

\newcommand{\mysection}[1]{\section{#1}\setcounter{equation}{0}}

\title{On upper bounds on expectations of gOSs based on DFR and DFRA distributions}
\author{Agnieszka Goroncy\\
Nicolaus Copernicus University\\
Chopina 12/18, 87100 Toru\'n, Poland\\
e-mail: gemini@mat.umk.pl
}

\date{}
\begin{document}

\maketitle{}

\begin{abstract}
We focus on the problem of establishing the optimal upper bounds on generalized order statistics which are based on the underlying cdf belonging to the family of distributions with decreasing failure rate and decreasing failure rate on the average. This issue has been previously considered by Bieniek [Projection bounds on expectations of generalized order statistics from DFR and DFRA families, {\it Statistics}, 2006; 40: 339--351], who established upper nonnegative mean-variance bounds with use of the projections of the compositions of density functions of the uniform generalized order statistic and the exponential distribution function onto the properly chosen convex cones. In this paper we obtain possibly negative upper bounds, by improving the zero bounds obtained by Bieniek for some particular cases of gOSs. We express the bounds in the scale units generated by the central absolute moments of arbitrary orders. We also describe the attainability conditions.
\vspace{2ex}

\noindent 2010 \textit{Mathematics Subject Classification}: 60E15,
62G32.

\vspace{2ex}

\noindent \textit{Key words: bound, generalized order statistics, decreasing failure rate, decreasing failure rate on the average, stochastic orderings}.
\end{abstract}

\mysection{Introduction}
Consider a random sample $X_1,\ldots,X_n$ of i.i.d. random variables with common cdf $F$ and finite absolute moment of order $p\in [1,\infty)$,
\begin{equation}
\sigma_p^p=\mathbb{E}|X_1-\mu|^p=\int\limits_0^1|F^{-1}(x)-\mu|^p dx,\label{pty_moment}
\end{equation}
with
\begin{equation}
\mu=\mathbb{E}X_1=\int\limits_0^1F^{-1}(x)dx.\label{mu}
\end{equation}
Let $U$ denote the standard uniform distribution function and $V$ the standard exponential distribution function, i.e. $V(x)=1-{\rm e}^{-x}$ with the density $v(x)={\rm e}^{-x}$, for $x\geq 0$.

We focus on the generalized order statistics $X_\gamma^{(1)},\ldots,X_\gamma^{(n)}$, introduced by Kamps (1995a, 1995b). They can be defined for an arbitrary vector of positive coefficients $\gamma=(\gamma_1,\ldots,\gamma_n)$ by the following quantile transformation
\begin{equation}\label{gOS_repr}
X^{(r)}_\gamma=F^{-1}\left(1-\prod\limits_{i=1}^r B_i\right),\quad 1\leq r\leq n,
\end{equation}
where $B_1,\ldots,B_n$ denote independent beta distributed random variables with distributions ${\rm Beta}(\gamma_1,1),\dots,{\rm Beta}(\gamma_n,1)$, respectively (see Cramer and Kamps, (2003)). The model of generalized order statistics (gOSs, for short) constitutes a unified approach containing many other popular models of ordered random variables, e.g. ordinary order statistics, record values, sequential order statistics, progressively censored type II order statistics and many more. \\
Let $U_\gamma^{(1)},\ldots,U_\gamma^{(n)}$ denote uniform gOSs, i.e. those arising from the standard uniform probability distribution function $U$. By $f_{\gamma,r}$ we denote the probability density function of $U^{(r)}_\gamma$, $1\leq r\leq n$, which, in the general case is given by
\begin{center}
$f_{\gamma,r}(x)=c_{r-1}G_r(x|\gamma_1,\ldots,\gamma_r)=c_{r-1}G_{r,r}^{r,0}$$\left(
\begin{tabular}{c|c}
    \multirow{2}{*}{$1-x$}&$\gamma_1,\ldots,\gamma_r$ \\
    & $\gamma_1-1,\ldots,\gamma_r-1$\\
\end{tabular}
\right)$, $\quad 0<x<1$,
\end{center}
where $c_{r-1}=\prod\limits_{j=1}^r\gamma_j$ and $G_{r,r}^{r,0}$ stands for the particular Meijer's G-function (see Mathai, Chapters 2 and 3, (1993)). We may further assume, that $\gamma_1\geq\ldots\geq\gamma_r>0$, since the marginal distribution of the single generalized order statistic does not depend on the ordering of the parameters $\gamma_1,\ldots,\gamma_r$.\\

This paper is devoted to the bounds on expectations of properly standardized gOSs, and its aim is to establish upper nonpositive (possibly negative) bounds on
\begin{equation}\label{problem}
\mathbb{E}\frac{X^{(r)}_\gamma-\mu}{\sigma_p},\quad p\geq 1,
\end{equation}
for the restricted families of the underlying distributions $F$, i.e. ones with decreasing failure rate or decreasing failure rate on the average. Such families of distributions can be defined in terms of the convex and star orders. First, we say that the distribution function $F$ succeeds some other fixed distribution function $W$ in the convex transform order, and write $F\succ_c W$, if the composition $F^{-1}\circ W$ is convex on the support of $W$. In particular, if $W=U$, then we say that $F$ belongs to the family of distributions with decreasing density (DD) and if $W=V$, then $F$ is an element of decreasing failure rate class of distributions (DFR). Indeed, since $F^{-1}(1-e^{-x})$ is convex on $[0,\infty)$, then the cumulative hazard function
$$
\Lambda(x)=V^{-1}(F(x))=-\ln[\bar{F}(x)],
$$
is concave on the support of $F$. Hence its derivative
$$
\lambda_F(x)=\Lambda'(x)=\frac{f(x)}{1-F(x)},
$$
called the failure rate function is nonincreasing. On the other hand, it is said that the distribution function $F$ belongs to the family with decreasing failure rate on the average distributions (DFRA) if $F\succ_* V$, i.e. $(F^{-1}(V(x))-F^{-1}(0))/x$ is nondecreasing on $[0,\infty)$.\\

There is a voluminous literature concerning the bounds on varieties of gOSs, in particular order statistics and their linear combinations (so called $L$-statistics), records and $k$th record values. The classical results on the bounds on the expected sample range expressed in terms of standard deviation units were determined by Plackett (1947). Later, at the beginning of the 1950s, Moriguti (1953) introduced the greatest convex minorant method which allowed to establish bounds for the arbitrary order statistics. A year later Gumbel (1954) and Hartley and David (1954) independently derived the generalized results on the expectation of the sample maximum and of the range, which were published in the same journal issue. Undoubtedly a turning point in the research was the introduction of the method of projection, proposed by Gajek and Rychlik (1996). It was immediately applied by Gajek and Rychlik (1998) in order to establish optimal bounds for the expectations of order statistics based on the life distributions with decreasing density or failure rate. Many applications of this method are comprehensively described in Rychlik (2001). Further, Rychlik (2002) considered optimal mean-variance bounds on order statistics from families of distributions determined by the star ordering, in particular distributions with decreasing density and failure rate on the average. Danielak (2003) obtained the bounds for expectations of the trimmed means from distributions with decreasing density and decreasing failure rate.

The general bounds on the expectations of gOSs can be found in Cramer et al. (2002) and in Goroncy (2014). Projection mean-variance bounds for gOSs from restricted families were established by Bieniek (2006, 2008) who considered distributions with decreasing failure rate and decreasing failure rate on the average, as well as distributions with decreasing density and decreasing density on the average. Recently Goroncy (2017) completed his results of the latter paper with the possibly negative upper bounds for some particular cases of gOSs based on the DD and DDA distributions.\\

\mysection{Main results}
We first focus on the case of distribution functions with decreasing failure rate, i.e., ones for which $F^{-1}\circ V$ are convex on $(0,\infty)$.
Note that we can express the expectations \re{problem} of the gOSs in terms of the following integral
\begin{equation}\label{EX_calka}
\mathbb{E}\frac{X^{(r)}_\gamma-\mu}{\sigma_p}=\int\limits_0^\infty\frac{F^{-1}(V(x))-\mu}{\sigma_p}(\hat{f}_{\gamma,r}(x)-1)v(x)dx,
\end{equation}
for $1\leq r\leq n$, and
\begin{equation}\label{f_tilde}
\hat{f}_{\gamma,r}(x)=f_{\gamma,r}(V(x)), \quad 0<x<\infty.
\end{equation}

Below we recall the results of Bieniek (see Theorem 4.2, (2006)), who established upper nonnegative bounds on the expectations of generalized order statistics based on DFR distributions. He applied the projection of the composition of the uniform gOS density with the exponential distribution function \re{f_tilde} onto a properly chosen convex cone. For this purpose we introduce
\begin{equation}\nonumber
\varrho_{j,r}=\sum\limits_{i=j}^r\frac{1}{\gamma_i},\quad 1\leq j\leq r,
\end{equation}
and the following function
\begin{equation}\nonumber
\alpha_*(y)=\frac{1}{2}\left(\sum\limits_{j=1}^{r-1}\frac{\varrho_{j,r}}{\gamma_j}\hat{f}_{\gamma,j}(y)+\left(\frac{1}{\gamma_r^2}-1\right)\hat{f}_{\gamma,r}(y)\right),\quad y\geq 0,
\end{equation}
which corresponds to the slope of the best linear approximation of the form $\hat{f}_{\gamma,r}(y)+\alpha(x-y)$ of the projected function \re{f_tilde} on $[y,\infty)$.
\begin{theorem}\label{th_Bieniek} If $r=1$ and $\gamma_1\geq 1$ then $EX_{\gamma}^{(1)}\leq \mu$. \\
Fix now $r\geq 2$ and parameters $\gamma_i\geq 1$, $i=1,\ldots,r$, where the second smallest of them is strictly greater than 1. Let $X_\gamma^{(r)}$ be the $r$th gOS based on the distribution function $F$ from DFR family of distributions with finite mean \re{mu} and positive variance \re{pty_moment} for $p=2$. \\
If $\varrho_{1,r}\leq 1$, then $EX_\gamma^{(r)}\leq \mu$.\\
If $1<\varrho_{1,r}\leq 2$, then
\begin{equation}\nonumber
\mathbb{E}\frac{X_\gamma^{(r)}-\mu}{\sigma_2}\leq \varrho_{1,r}-1,
\end{equation}
and the bound is attained for the following exponential distribution
\begin{equation}\nonumber
F(x)=\left\{
\begin{array}{ll}
0,& {\rm if}\hspace{3mm} x\leq \mu-\sigma_2,\\
&\\
1-{\rm exp}\left(-\dfrac{x-\mu}{\sigma_2}-1\right),& {\rm if}\hspace{3mm} x>\mu-\sigma_2.
\end{array}
\right.
\end{equation}
If $\varrho_{1,r}>2$, then
\begin{equation}\label{osz_Bieniek}
\mathbb{E}\frac{X_{\gamma}^{(r)}-\mu}{\sigma_2}\leq C=C_{1,r}(\gamma),
\end{equation}
where
\begin{equation}\nonumber
C^2=\int\limits_{0}^{y^*}(\hat{f}_{\gamma,r}(x))^2{\rm e}^{-x}dx+{\rm e}^{-y^*}\left\{(\hat{f}_{\gamma,r}(y^*))^2+2\alpha_*(y^*)\hat{f}_{\gamma,r}(y^*)+2(\alpha_*(y^*))^2\right\}-1,
\end{equation}
where $y^*$ is the smallest positive solution to the following equation
\begin{equation}\nonumber
\sum\limits_{j=1}^{r-1}\frac{1}{\gamma_j}\left(1-\frac{\varrho_{j,r}}{2}\right)\hat{f}_{\gamma,j}(y)-\frac{(\gamma_r-1)^2}{2\gamma_r^2}\hat{f}_{\gamma,r}(y)=0.
\end{equation}
The bound \re{osz_Bieniek} is attained for the distribution function
\begin{equation}\nonumber
F(x)=\left\{
\begin{array}{ll}
0,& {\rm if}\hspace{3mm} \dfrac{x-\mu}{\sigma_2}\leq -\dfrac{1}{C},\\
&\\
f_{\gamma,r}^{-1}\left(C\dfrac{x-\mu}{\sigma_2}+1\right),& {\rm if}\hspace{3mm} -\dfrac{1}{C}<\dfrac{x-\mu}{\sigma_2}\leq\dfrac{\hat{f}_{\gamma,r}(y^*)-1}{C},\\
&\\
V\left(C\dfrac{x-\mu}{\sigma_2\alpha_*(y^*)}+\dfrac{1-\hat{f}_{\gamma,r}(y^*)}{\alpha_*(y^*)}+y^*\right), & {\rm if}\hspace{3mm} \dfrac{x-\mu}{\sigma_2}>\dfrac{\hat{f}_{\gamma,r}(y^*)-1}{C}.
\end{array}
\right.
\end{equation}
\end{theorem}
Note that for cases $r\geq 1$ with $\varrho_{1,r}\leq 1$, Bieniek derived zero bounds without describing the attainability conditions. In general, the positivity and negativity of the upper bounds on the expected gOSs depends on the parameters $\gamma_1,\ldots,\gamma_r$, and the restrictions imposed on the parent distribution function $F$. The projection method is appropriate in cases when it results with nonconstant projections (and positive upper bounds) of functions \re{f_tilde}, which depend basically on their shapes. Otherwise this procedure returns zero bounds and another method should be employed in order to obtain possibly negative ones. Our objective is to improve Bieniek's bounds in some particular cases.

In the first proposition we present two most general cases of optimal nonpositive bounds, which can either be zero or strictly negative.
\begin{proposition}
Let $\mathbb{E}X^{(r)}_{\gamma}$ be the $r$th generalized order statistic based on parameter vector $\gamma=(\gamma_1\ldots,\gamma_r)\in \mathbb{R}_+^r$, and a DFR parent distribution function $F$ with expectation $\mu$.
\begin{itemize}
\item[(i)] If $\varrho_{1,r}=\sum\limits_{i=1}^r\frac{1}{\gamma_i}=1$, then the bound $\mathbb{E}X^{(r)}_{\gamma}\leq \mu$ is attained by the exponential baseline distribution functions.
\item[(ii)] If $0<\varrho_{1,r}<1$, then for every DFR baseline distribution function $F$ we have $\mathbb{E}X^{(r)}_{\gamma}<\mu$.
\end{itemize}
\end{proposition}
\textsc{Proof.} Formula \re{gOS_repr} can be rewritten as
$$
X^{(r)}_\gamma=F^{-1}\left(1-\prod\limits_{i=1}^r (1-U_i)^{1/\gamma_i}\right),\quad 1\leq r\leq n,
$$
where $U_i$, $i=1,\ldots,r$, are i.i.d. standard uniform random variables. Further on we can represent them as $U_i=1-{\rm exp}(-V_i)$, $i=1,\ldots,r$, with $V_i$ being i.i.d. standard exponential. Therefore
$$
X^{(r)}_\gamma=F^{-1}\left(1-\prod\limits_{i=1}^r {\rm exp}\left(\frac{V_i}{\gamma_i}\right)\right)=F^{-1}\left(V\left(\sum\limits_{i=1}^r\frac{V_i}{\gamma_i}\right)\right).
$$
If $F(x)=V(\frac{x-\theta}{\lambda})$ is the exponential distribution function with location parameter $\theta$ and scale parameter $\lambda>0$, and $\sum\limits_{i=1}^r\frac{1}{\gamma_i}=1$, then
$$
\mathbb{E}X^{(r)}_{\gamma}=\sum\limits_{i=1}^r\frac{\mathbb{E}(\lambda V_i+\theta)}{\gamma_i}=\lambda+\theta=\mu.
$$
Suppose now that $0<\sum\limits_{i=1}^r\frac{1}{\gamma_i}=c<1$. Then obviously $\sum\limits_{i=1}^r\frac{1}{c\gamma_i}=1$ and $\sum\limits_{i=1}^r\frac{V_i}{\gamma_i}<\sum\limits_{i=1}^r\frac{V_i}{c\gamma_i}$. Composition $F^{-1}\circ V$ for $F$ being DFR is possibly first constant, and then strictly increasing. Consequently,
$$
F^{-1}\left(V\left(\sum\limits_{i=1}^r\frac{V_i}{\gamma_i}\right)\right)\leq F^{-1}\left(V\left(\sum\limits_{i=1}^r\frac{V_i}{c\gamma_i}\right)\right),
$$
and the inequality is strict for sufficiently large arguments. Since any positive combination  of $V_i$ has a positive probability measure of these arguments, it follows that
$$
\mathbb{E}X^{(r)}_{\gamma}=\mathbb{E}F^{-1}\left(V\left(\sum\limits_{i=1}^r\frac{V_i}{\gamma_i}\right)\right)<\mathbb{E}F^{-1}\left(V\left(\sum\limits_{i=1}^r\frac{V_i}{c\gamma_i}\right)\right)=\mathbb{E}X^{(r)}_{c\gamma}\leq \mu,
$$
where the last inequality follows from Theorem \ref{th_Bieniek}.\hfill$\Box$\\

In the next proposition we present the most general results for the upper negative bounds on \re{problem} with $1\leq p<\infty$.
\begin{proposition}\label{dfr_prop1}
Fix $r\geq 1$ and $\gamma_i>1$, $i=1,\ldots,r$, are such that $\varrho_{1,r}=\sum\limits_{i=1}^r\frac{1}{\gamma_i}< 1$. Let $X_\gamma^{(r)}$ be the $r$th gOS based on the distribution function $F$ from DFR family of distributions with finite mean \re{mu} and absolute central moment \re{pty_moment} with $1\leq p<\infty$. We then have the following bound
\begin{equation}\label{bound_p}
\mathbb{E}\frac{X_\gamma^{(r)}-\mu}{\sigma_p}\leq -\inf\limits_{0\leq\alpha<\infty}B_p(\alpha),
\end{equation}
where
\small{
\begin{equation}\label{B_p}
B_p(\alpha)=\dfrac{b(\alpha)^{-1}}{\left[{\rm e}^{-\alpha p}(1-{\rm e}^{-\alpha})+\int\limits_\alpha^{\alpha+{\rm e}^{-\alpha}}(\alpha+{\rm e}^{-\alpha}-x)^p{\rm e}^{-x}dx+\int\limits_{\alpha+{\rm e}^{-\alpha}}^\infty(x-{\rm e}^{-\alpha}-\alpha)^p{\rm e}^{-x}dx\right]^{1/p}},
\end{equation}
}
with
\begin{equation}\label{b_alpha}
b(\alpha)={\rm e}^{\alpha}\left(1-\sum\limits_{j=1}^r\dfrac{\varrho_{j,r}}{\gamma_j}\hat{f}_{\gamma,j}(\alpha)\right)^{-1}.
\end{equation}
The equality holds (possibly is approached in the limit) for the mixtures of atoms at $c_2$ and shifted exponential distributions with the following cdf
\begin{equation}\label{F_prop1}
F_\alpha(x)=\left\{
\begin{array}{ll}
0,&x<c_2,\\
1-{\rm exp}\left(-\dfrac{x-c_2}{c_1}-\alpha\right),& x\geq c_2,
\end{array}
\right.
\end{equation}
with probabilites $1-{\rm e}^{-\alpha_0}$ and ${\rm e}^{-\alpha_0}$, respectively, for
\begin{eqnarray}
c_1&=&b(\alpha_0)B_p(\alpha_0)\sigma_p,\nonumber\\
c_2&=&\mu-{\rm e}^{-\alpha_0}c_1,\nonumber
\end{eqnarray}
where $\alpha_0$ is the parameter in which the infimum of \re{B_p} is attained (possibly approached in the limit).

\end{proposition}
\textsc{Proof.} The proof is analogous to the proofs presented in papers of Goroncy (2017, Proposition 1) and Rychlik (2009, Theorem 2), therefore we present only the main steps of the reasoning. This particular case requires consideration of the composition of the density function of the uniform $r$th gOS and the exponential distribution function $V$, which is given by \re{f_tilde}, and $1\leq p<\infty$. \\
Let us rewrite \re{EX_calka} as the following functional
\begin{equation}\nonumber
T_h(g)=\int\limits_0^\infty h(x)g(x)v(x)dx,
\end{equation}
represented by
$$
h(x)=\hat{f}_{\gamma,r}(x)-1,
$$
with
\begin{eqnarray}\label{g}
g(x)=\frac{F^{-1}(V(x))-\mu}{\sigma_p}.
\end{eqnarray}
Further, consider a convex and bounded subset $\mathcal{G}_p$ of $L^p([0,\infty),v(x)dx)$, which consists of nondecreasing and convex functions that integrate to zero and have the unit $p$th norm. Obviously, since $F$ belongs to DFR family of distributions, functions \re{g} are elements of $\mathcal{G}_p$. According to Bieniek (2006, Theorem 4.2) and our setting of parameters, we assume that $T_h(g) < 0$. Transforming $g$ into $\tilde{g}=-g/T_h(g)$ we also transform set $\mathcal{G}_p$ into $\tilde{\mathcal{G}}_p$, which is a set that consists of all the nondecreasing and convex functions $\tilde{g}$ of $L^p([0,\infty),v(x)dx)$, which satisfy
\begin{eqnarray}
T_1(\tilde{g})&=&\int\limits_0^\infty \tilde{g}(x)v(x)dx=0,\label{T1}\\
T_h(\tilde{g})&=&\int\limits_0^\infty h(x)\tilde{g}(x)v(x)dx= -1.\label{T2}
\end{eqnarray}
Noticing that $||\tilde{g}||_p=-1/T_h(g)$, for $\tilde{g}\in\tilde{\mathcal{G}}_p$ and $g\in\mathcal{G}_p$, which increases as $T_h(g)$ increases, we conclude that the problem of finding the supremum of $T_h$ over $\mathcal{G}_p$ and the problem of determining
\begin{equation}\label{norma_tilde_g}
\sup\{||\tilde{g}||_p, \tilde{g}\in \mathcal{\tilde{G}}_p\},
\end{equation}
are equivalent. Since the norm is a convex functional, we may simplify the calculations by confining ourselves only to the extreme elements of $\tilde{\mathcal{G}}_p$, i.e. maximizing the norm \re{norma_tilde_g} on the subset of functions of the form
\begin{equation}\label{g_alpha}
\tilde{g}_\alpha(x)=a(\alpha)+b(\alpha)(x-\alpha)\mathbf{1}_{[\alpha,\infty)}(x), \quad \alpha\in(0,\infty),
\end{equation}
with coefficients $a(\alpha)$ and $b(\alpha)$ complying with conditions \re{T1} and \re{T2}. The substantiation of the fact that extreme elements of $\tilde{\mathcal{G}}_p$ are broken lines defined in \re{g_alpha}, is exactly the same as in Rychlik (2009, pp.59--60). Therefore we claim that
\begin{equation}\label{supremum_E}
\sup\mathbb{E}\frac{X^{(r)}_\gamma-\mu}{\sigma_p}=\sup\{T_h(g), g\in \mathcal{G}_p\}=-\inf\limits_{0\leq\alpha<\infty}\frac{1}{||\tilde{g}_\alpha||_p}.
\end{equation}
We first calculate the coefficients of \re{g_alpha}. Using the following integral calculations obtained by Bieniek (see Lemma 3.1, (2006))
\begin{eqnarray}
\int\limits_\alpha^\infty\hat{f}_{\gamma,r}(x){\rm e}^{-x}dx&=&{\rm e}^{-\alpha}\sum\limits_{j=1}^r\frac{1}{\gamma_j}\hat{f}_{\gamma,j}(\alpha),\nonumber\\
\int\limits_\alpha^\infty(x-\alpha)\hat{f}_{\gamma,r}(x){\rm e}^{-x}dx&=&{\rm e}^{-\alpha}\sum\limits_{j=1}^r\frac{\varrho_{j,r}}{\gamma_j}\hat{f}_{\gamma,j}(\alpha),\nonumber
\end{eqnarray}
we obtain that
\begin{eqnarray}
T_1(\tilde{g}(\alpha))&=&\int\limits_0^\infty\tilde{g}_\alpha(x){\rm e}^{-x}dx=a(\alpha)+b(\alpha){\rm e}^{-\alpha},\\
T_h(\tilde{g}(\alpha))&=&\int\limits_0^\infty(\hat{f}_{\gamma,r}(x)-1)\tilde{g}_\alpha(x){\rm e}^{-x}dx=b(\alpha){\rm e}^{-\alpha}\left(\sum\limits_{j=1}^r\frac{\varrho_{j,r}}{\gamma_j}\hat{f}_{\gamma,j}(\alpha)-1\right),
\end{eqnarray}
and together with conditions \re{T1} and \re{T2}, conclude that
\begin{equation}\nonumber
a(\alpha)=-b(\alpha){\rm e}^{-\alpha},
\end{equation}
where $b(\alpha)$ is given by \re{b_alpha}.
We are now in the position of determining the norm of \re{g_alpha}, which is required for establishing \re{supremum_E},
\begin{eqnarray}
||\tilde{g}_\alpha||_p^p&=&\int\limits_0^\infty\left|a(\alpha)+b(\alpha)(x-\alpha)\mathbf{1}_{[\alpha,\infty)}(x)\right|^p{\rm e}^{-x}dx= b^p(\alpha)\left[{\rm e}^{-p\alpha}(1-{\rm e}^{-\alpha})\right.\nonumber\\
&& +\int\limits_\alpha^{\alpha+{\rm e}^{-\alpha}}(\alpha+{\rm e}^{-\alpha}-x)^p{\rm e}^{-x}dx
+\int\limits_{\alpha+{\rm e}^{-\alpha}}^\infty(x-\alpha-{\rm e}^{-\alpha})^p{\rm e}^{-x}dx],\nonumber
\end{eqnarray}
and eventually results in the bound \re{B_p}.

The equality in \re{bound_p} is attained for the distribution functions $F_\alpha$, which satisfy the following condition
$$
\frac{F_\alpha^{-1}(V(x))-\mu}{\sigma_p}=\frac{\tilde{g}_\alpha(x)}{||\tilde{g}_\alpha||_p}.
$$
The above proportion can be rewritten as
$$
F^{-1}_\alpha(1-{\rm e}^{-x})=\left\{
\begin{array}{ll}
\mu-\sigma_pb(\alpha){\rm e}^{-\alpha}B_p(\alpha),&0\leq x<\alpha,\\
&\\
\mu+\sigma_pb(\alpha)B_p(\alpha)(x-\alpha-{\rm e}^{-\alpha}),&\alpha\leq x<\infty,
\end{array}
\right.
$$
which is consistent with the equality conditions described in the proposition.
\hfill$\Box$\\

Note that it is not easy to determine the integrals in \re{B_p} of Proposition \ref{dfr_prop1} in general, especially for noninteger $p$. However, we managed to confirm the zero bound and specify the equality conditions for the special case of the first generalized order statistic, i.e. $r=1$ and arbitrarily chosen $1<p<\infty$. This is possible due to the fact that the density of gOS in this case can be easily determined and that simplifies calculations. These results are presented in the proposition below, while the special case $r=1$ with $p=1$ is considered further in a separate corollary.
\begin{proposition}\label{Prop_2}
Let $X_\gamma^{(1)}$ be the first generalized order statistic with parameter $\gamma=\gamma_1\geq1$, based on the distribution function $F$ from DFR family of distributions with mean $\mu\in \mathbb{R}$ and finite $p$th central absolute moment $\sigma_p$, for $1<p<\infty$. If $\gamma=1$, then $\mathbb{E}X^{(1)}_1=\mu$, which holds for any DFR distribution function $F$. If $\gamma>1$, then we have the following upper bound
$$
\mathbb{E}\frac{X^{(1)}_\gamma-\mu}{\sigma_p}\leq 0.
$$
The equality above is attained in limit by sequences of distributions $F_\alpha$, $0<\alpha<1$, which are the mixtures of atoms in $\mu-\frac{\sigma_p\alpha}{N_p(\alpha)}$ and exponential distributions with the location $\mu-\frac{\sigma_p\alpha}{N_p(\alpha)}$ and scale $\frac{\sigma_p}{N_p(\alpha)}$ parameters with probabilities $1-\alpha$ and $\alpha$ respectively, for $\alpha\longrightarrow 0$, and
$$
N_p(\alpha)=\left[\alpha^p-\alpha^{p+1}+\alpha{\rm e}^{-\alpha}\left(\int\limits_{0}^\alpha y^p{\rm e}^ydy+\Gamma(p+1)\right)\right]^{1/p}.
$$
\end{proposition}

\textsc{Proof.} The case $r=1$ and $\gamma=1$, is obvious, since $\mathbb{E}X^{(1)}_{1}=\mathbb{E}X_1=\mu$ for any parent distribution $F$.

Suppose $\gamma>1$. Rychlik (see proof of Theorem 1, (2009)), who obtained analogous results for the common order statistics, has shown that distribution functions $F_\alpha$, $0<\alpha<1$, which are obviously DFR, have mean equal to \re{mu} and $p$th absolute central moment equal to \re{pty_moment}. Moreover, the quantile function of such distributions is given by
$$
F_\alpha^{-1}(x)=\left\{
\begin{array}{ll}
\mu-\frac{\sigma_p}{N_p(\alpha)}\alpha,&0<x<1-\alpha,\\
&\\
\mu+\frac{\sigma_p}{N_p(\alpha)}\left(-\ln(1-x)+\ln\alpha-\alpha\right),&1-\alpha\leq x<1.
\end{array}
\right.
$$
 It is sufficient to check that
\begin{equation}\nonumber
\mathbb{E}\frac{X^{(1)}_\gamma-\mu}{\sigma_p} \longrightarrow 0, \quad {\rm as}\hspace{2mm} \alpha\longrightarrow0.
\end{equation}
First note that for the first uniform gOS we have the following density and the cumulative probability functions
\begin{eqnarray}
f_{\gamma,1}(x)=\gamma_1(1-x)^{\gamma_1-1},\nonumber\\
F_{\gamma,1}(x)=1-(1-x)^{\gamma_1},\nonumber
\end{eqnarray}
respectively, where $0\leq x<1$. Therefore using the formula
\begin{equation}\nonumber
\int\limits_\alpha^\infty(x-\alpha)\hat{f}_{\gamma,1}(x){\rm e}^{-x}dx=\frac{1}{\gamma_1}{\rm e}^{-\alpha\gamma_1},
\end{equation}
we obtain
\begin{eqnarray}
0&>&\mathbb{E}_{F_\alpha}\frac{X_\gamma^{(1)}-\mu}{\sigma_p}=\int\limits_0^1\frac{F_\alpha^{-1}(x)-\mu}{\sigma_p}(f_{\gamma,1}(x)-1)dx\nonumber\\
&&=\frac{1}{N_p(\alpha)}\int\limits_{1-\alpha}^1(-\ln(1-x)+\ln\alpha)(f_{\gamma,1}(x)-1)dx\nonumber\\
&&=\frac{1}{N_p(\alpha)}\int\limits_{-\ln\alpha}^\infty(y+\ln\alpha)(f_{\gamma,1}(1-{\rm e}^{-y})-1){\rm e}^{-y}dy\nonumber\\
&&=\frac{1}{N_p(\alpha)}\left(\int\limits_{-\ln\alpha}^\infty(y+\ln\alpha)f_{\gamma,1}(1-{\rm e}^{-y}){\rm e}^{-y}dy-\alpha\int_0^\infty y{\rm e}^{-y}dy\right)\nonumber\\
&&=\alpha^{1-\frac{1}{p}}\dfrac{\frac{1}{\gamma_1}\alpha^{\gamma_1-1}-1}{\left\{\alpha^{p-1}-\alpha^{p}+{\rm e}^{-\alpha}\left[\int\limits_0^\alpha y^p{\rm e}^{y}dy+\Gamma(p+1)\right] \right\}^{1/p}},\nonumber
\end{eqnarray}
which tends to 0 if $\alpha\longrightarrow 0$, and completes the proof. \hfill$\Box$\\

Note that the first generalized order statistic has the distribution identical with the minimum of $n$ i.i.d. random variables with the distribution function $1-(1-F(x))^{\gamma_1/n}$, which shares the DFR property. Therefore, knowing the bounds for the ordinary order statistics, we also obtain the bounds for the first generalized order statistics, with use of this transformation. This is why results of Proposition \ref{Prop_2} above coincide with the result of Rychlik (see Theorem 1, (2009)), who considered upper nonpositive bounds on low rank common order statistics in the same model of DFR family of distributions.

The case of bounds given in scale units generated by the mean absolute deviation $\sigma_1$, is treated separately in the corollary below, which follows from Proposition \ref{dfr_prop1}. Note that in this case
\begin{equation}\label{g_1}
||\tilde{g}_\alpha||_1=2{\rm e}^{-\alpha-{\rm e}^{-\alpha}}b(\alpha),
\end{equation}
since the denominator of \re{B_p} reduces to $2{\rm e}^{-\alpha-{\rm e}^{-\alpha}}$.

\begin{corollary}\label{tw_p1}
Fix $r\geq 1$ and $\gamma_i>1$, $i=1,\ldots,r$, are such that $\varrho_{1,r}=\sum\limits_{i=1}^r\frac{1}{\gamma_i}< 1$. Let $X_\gamma^{(r)}$ be the $r$th gOS based on the distribution function $F$ from DFR family of distributions with finite absolute mean deviation $\sigma_1$ defined in \re{pty_moment} for $p=1$. We have the following bound
\begin{equation}\nonumber
\mathbb{E}\frac{X^{(r)}_\gamma-\mu}{\sigma_1}\leq -\inf\limits_{0\leq\alpha<\infty}B_1(\alpha),
\end{equation}
where
\begin{equation}\label{B_1}
B_1(\alpha)=\frac{1}{2}{\rm e}^{{\rm e}^{-\alpha}}\left[1-\sum\limits_{j=1}^r\frac{\varrho_{j,r}}{\gamma_j}\hat{f}_{\gamma,j}(\alpha)\right],\quad 0\leq\alpha<\infty.
\end{equation}
The equality holds (possibly is approached in the limit) for mixtures of atoms in $\mu-\frac{1}{2}{\rm e}^{{\rm e}^{-\alpha_0}}\sigma_1$ and shifted exponential distributions \re{F_prop1} with $c_1=\frac{1}{2}{\rm e}^{{\rm e}^{-\alpha_0}+\alpha_0}\sigma_1$ and $c_2=\mu-\frac{1}{2}{\rm e}^{{\rm e}^{-\alpha_0}}\sigma_1$
with probabilites $1-{\rm e}^{-\alpha_0}$ and ${\rm e}^{-\alpha_0}$, respectively, where $\alpha_0$ is the parameter in which the infimum of \re{B_1} is attained (possibly approached in the limit).
\end{corollary}

In particular, bounds for the first generalized order statistics, given in the absolute mean deviation units are strictly negative, presented in the corollary below.

\begin{corollary}
Let $X_\gamma^{(1)}$ be the first generalized order statistic with parameter $\gamma=\gamma_1> 1$, based on the distribution function $F$ from DFR family of distributions with  finite mean absolute deviation $\sigma_1$. We have the following upper bound
\begin{equation}\label{dfr_p11}
\mathbb{E}\frac{X^{(1)}_\gamma-\mu}{\sigma_1}\leq -\inf\limits_{0<\beta\leq 1}\frac{1}{2}{\rm e}^{\beta}\left(1-\frac{1}{\gamma}\beta^{\gamma-1}\right).
\end{equation}
The equality conditions in \re{dfr_p11} correspond to those described in Corollary \ref{tw_p1}, with $\alpha_0=-\ln\beta_0$, where $\beta_0$ is the point in which the infimum of the rhs of \re{dfr_p11} is attained.\end{corollary}

Note that here we have
\begin{equation}\nonumber
b(\alpha)={\rm e}^{\alpha}\left(1-\frac{1}{\gamma_1}{\rm e}^{-\alpha(\gamma_1-1)}\right)^{-1},
\end{equation}
which for $0\leq\beta={\rm e}^{-\alpha}<1$, together with \re{g_1} results in the desired bound \re{dfr_p11}.\\

The numerical values of bounds for the first gOSs expressed in $\sigma_1$ units for some particular cases of the parameter $\gamma_1$ are presented in Table 1. Parameter $\beta_0$ stands for the argument of the rhs of \re{dfr_p11}, for which the infimum is attained.
\begin{center}
\begin{table}[ht]
  \begin{tabular}{|l|l|l||l|l|l|}
    $\gamma_1$ & $\beta_0$& bound &$\gamma_1$ & $\beta_0$&bound\\\hline
    1,005 &0,99& -0,0068 &1,2&0,6461&-0,2255\\
    1,01 &0,9801&  -0,0135&1,3&0,4980&-0,3093\\
     1,03&0,9412&-0,0396&1,4&0,3661&-0,3765\\
1,04&0,9221&-0,0523  &1,5&0,2500&-0,4280\\
1,05&0,9032&-0,0647&1,6&0,1509&-0,4646\\
1,06&0,8846&-0,0769&1,7&0,0720&-0,4872\\
1,07&0,8662&-0,0889&1,8&0,0196&-0,4977\\
1,08&0,8480&-0,1006&1,9&0,0006&-0,4999\\
1,09&0,8301&-0,1122&2&0&-0,5000\\
1,1&0,8123&-0,1235&3&0&-0,5000\\
  \end{tabular}
  \caption{Bounds on expectations of the standardized $1^{st}$ gOSs, $X_{\gamma_1}^{(1)}$, given in the absolute mean deviation units $\sigma_1$, DFR case.}
\end{table}
\end{center}
The numerical simulations show that for increasing values of $1<\gamma_1<2$, the value of $\beta_0$ decreases from 1 to 0 while the bound gradually decreases to $-1/2$. It eventually stabilizes its value in $-1/2$ beginning with $\gamma_1=2$, when the infimum of the rhs of \re{dfr_p11} is attained in $\beta_0=0$. This coincides with results of Goroncy (2014, Theorem 4) in the general case of the arbitrary underlying distribution function for $p=1$.\\

Let us now consider gOSs which are based on DFRA distribution functions. The case of upper nonnegative bounds on \re{problem} for $p=2$ was solved by Bieniek (2006, Theorem 4.3). However, for some particular cases of gOSs, the projection method, which was used there, resulted in zero bounds. Recently, Goroncy (2017) presented bounds on expectations of the standardized gOS arising from decreasing density on the average distributions. This method of obtaining possibly negative bounds will be applied here in order to obtain respective results for decreasing failure rate on the average distributions.

The results presented below are just the immediate application of the Proposition 3 of Goroncy (2017) in case when a particular condition presented in Bieniek (2006, Theorem 4.3) is satisfied.

\begin{proposition}
Let $X_{\gamma}^{(r)}$, $r\geq 2$, be the $r$th gOS based on the distribution function $F$ which belongs to DFRA family of distributions and fix $1\leq p<\infty$.
Let $\hat{\beta}$ be the only zero of the equation $\hat{f}_{\gamma,r}(\beta)=1$ in the interval $(0,\hat{\theta})$, for $\hat{\theta}$ being the smallest inflection point of \re{f_tilde}. If the following condition is satisfied
$$
\sum\limits_{j=1}^r\frac{1}{\gamma_j}(\varrho_{j,r}+\hat{\beta})\hat{f}_{\gamma,r}(\hat{\beta})\leq 1+\hat{\beta},
$$
then the following bound holds
\begin{equation}\label{oszacowanie_DGFRA}
\mathbb{E}\frac{X^{(r)}_\gamma-\mu}{\sigma_p}\leq -\inf\limits_{0<\alpha<\infty}B^*_p(\alpha),
\end{equation}
where
\begin{eqnarray}
B^*_p(\alpha)&=&\left\{
\begin{array}{ll}
\frac{1}{b_\alpha}\left[((\alpha+1){\rm e}^{-\alpha})^p(1-{\rm e}^{-\alpha})+\int\limits_{\alpha}^{(\alpha+1){\rm e}^{-\alpha}}(-x+(\alpha+1){\rm e}^{-\alpha})^p{\rm e}^{-x}dx\right.&\\
\left.+\int\limits_{(\alpha+1){\rm e}^{-\alpha}}^\infty(x-(\alpha+1){\rm e}^{-\alpha})^p{\rm e}^{-x}dx\right]^{-1/p},&0<\alpha<\alpha_0,\\
&\\
\frac{1}{b_\alpha}\left[\right((\alpha+1){\rm e}^{-\alpha})^p(1-{\rm e}^{-\alpha})+\int\limits_{\alpha}^\infty(x-(\alpha+1){\rm e}^{-\alpha})^p{\rm e}^{-x}dx]^{-1/p},&\alpha\geq\alpha_0,
\end{array}
\right.\nonumber
\\
b_\alpha&=&-{\rm e}^\alpha\left[\sum\limits_{j=1}^r\frac{1}{\gamma_j}\hat{f}_{\gamma,j}(\alpha)(\alpha+\varrho_{j,r})-\alpha-1\right]^{-1},\nonumber
\end{eqnarray}
with $\alpha_0\backsimeq0,8065$ being the only solution of the equation $\alpha=(\alpha+1){\rm e}^{-\alpha}$ in $(0,\infty)$.
The equality in \re{oszacowanie_DGFRA} is attained for distribution functions $F_{\alpha_*}$ such that
\begin{equation}\label{rownosc_DGFRA}
F^{-1}_{\alpha_*}(V(x))=\left\{
\begin{array}{ll}
  -b_{\alpha_*}(\alpha_*+1){\rm e}^{-\alpha_*}B^*_p(\alpha_*)\sigma_p +\mu, & 0<x<\alpha_*, \\
  &\\
  (x-(\alpha_*+1){\rm e}^{-\alpha_*})b_{\alpha_*} B^*_p(\alpha_*)\sigma_p+\mu, & \alpha_*<x<d.
\end{array}
\right.
\end{equation}
where $\alpha_*$ is the argument for which the infimum of $B^*_p$ is attained. If the infimum of $B^*_p$ is attained in the limit, then the equality is also attained in the limit by sequences of distributions described in \re{rownosc_DGFRA}.
\end{proposition}
Numerical calculations show, that in some cases (e.g. for particular settings of progressively type II censored order statistics), the infimum of $B^*_p$ is approached in limit for $\alpha\rightarrow\infty$ and equals to zero. This is not surprising, but it is not easy to show analytically that in general the zero infimum of $B^*_p$ is attained in the limit, at the right exponential support interval.

\section*{Acknowledgements}
The research was supported by the Polish National Science Center Grant\break
no. 2011/01/D/ST1/04172. The author is very grateful to anonymous referees, whose helpful comments and remarks helped to significantly improve the manuscript, especially to the one who suggested addition of Proposition 1.


\begin{thebibliography}{30}
\bibitem{} Bieniek, M., (2006), Projection bounds on expectations of generalized order statistics from DFR and DFRA families, {\it Statistics}, \textbf{40}: 339--351.
\bibitem{} Bieniek, M., (2008), Projection bounds on expectations of generalized order statistics from DD and DDA families, {\it J. Statist. Plann. Inference},\textbf{ 138}: 971--981.
\bibitem{} Cramer, E., Kamps, U., (2003), Marginal distributions of sequential and generalized order statistics, \emph{Metrika}, \textbf{58}: 293--310.
\bibitem{} Cramer, E., Kamps, U., Rychlik, T., (2002), Evaluations of expected generalized order statistics in various scale units, \emph{Appl. Math.} (Warsaw) \textbf{29}: 285--295.
\bibitem{} Danielak, K., (2003), Sharp upper mean-variance bounds for trimmed means from restricted families, \textit{Statistics}, \textbf{37}: 305--324.
\bibitem{} Gajek, L., Rychlik, T., (1996), Projection method for moment bounds on order statistics from restricted families. I. Dependent case. \textit{J. Multivar. Anal.}, \textbf{57}: 156--174.
\bibitem{} Gajek, L., Rychlik, T., (1998), Projection method for moment bounds on order statistics from restricted families. II. Independent case, \textit{J. Multivar. Anal.}, \textbf{64}: 156--182.
\bibitem{} Goroncy, A., (2014), Bounds on expected generalized order statistics, \textit{Statistics}, \textbf{48}: 593--608.
\bibitem{} Goroncy, A., (2017), Upper non-positive bounds on expectations of generalized order statistics from DD and DDA populations, \textit{Comm. Stat. Theor. Meth.}, \textbf{46}:24: 11972--11987.
\bibitem{} Gumbel, E.J., (1954), The maxima of the mean largest value and of the range, \textit{Ann. Math. Stat.}, \textbf{25}: 76--84.
\bibitem{} Hartley, H.O., David, H.A., (1954), Universal bounds for mean range and extreme observation, \textit{Ann. Math. Stat.}, \textbf{25}: 85--99.
\bibitem{} Kamps, U., (1995a), \textit{A Concept of Generalized Order Statistics}, Teubner, Stuttgart.
\bibitem{} Kamps, U., (1995b), A concept of generalized order statistics, \textit{J. Statist. Plann. Inference}, \textbf{48}: 1--23.
\bibitem{} Mathai, A. M., (1993), \textit{A Handbook of Generalized Special Functions for Statistical and Physical Sciences}, Clarendon Press, New York.
\bibitem{} Moriguti, S., (1953), A modification of Schwarz's inequality with applications to distributions. {\it Ann. Math. Statist.} {\bf 24}: 107--113.
\bibitem{} Plackett, R.L., (1947), Limits of the ratio of mean range to standard deviation. {\it Biometrika} {\bf 34}: 120--122.
\bibitem{} Rychlik, T., (2001), \textit{Projecting statistical functionals}, Springer-Verlag, New York.
\bibitem{} Rychlik, T., (2002), Optimal mean-variance bounds on order statistics from families determined by star ordering, {\it Appl. Math.}, Warsaw, \textbf{29}: 15--32.
\bibitem{} Rychlik, T., (2009), Non-positive upper bounds on expectations of low rank order statistics from DFR populations, {\it Statistics}, \textbf{43}: 53--63.
\end{thebibliography}
\end{document}